\documentclass[12pt]{amsart}

\usepackage{enumitem}
\usepackage{psfrag}
\usepackage{graphicx}
\usepackage{pinlabel}
\usepackage{graphicx}
\usepackage{amsmath}
\usepackage{amssymb}
\usepackage{amscd}
\usepackage{autobreak}
\usepackage{picinpar}
\usepackage{times}
\usepackage{pb-diagram}
\usepackage{graphicx}
\usepackage{wrapfig}
\usepackage{xspace}
\usepackage{hyperref}
\usepackage{url}
\usepackage{caption}
\usepackage{subcaption}
\usepackage{fancyhdr}
\usepackage{overpic}
\usepackage{color}
\usepackage[square,comma,sort&compress,numbers]{natbib}
\usepackage{latexsym}
\usepackage{mathrsfs}
\usepackage{amsfonts}
\usepackage{hyperref}
\usepackage{amsthm}
\usepackage{appendix}
\usepackage{pdfsync}

\theoremstyle{definition}
\newtheorem{theorem}{Theorem}[section]
\newtheorem{lemma}[theorem]{Lemma}
\newtheorem{proposition}[theorem]{Proposition}

\newtheorem{definition}[theorem]{Definition}

\newtheorem{remark}[theorem]{Remark}
\newtheorem*{theorem*}{Theorem}
\def\qed{\hfill{Q.E.D.}\smallskip}

\begin{document}

\title{\bf Combinatorial curvature flows with surgery for inversive distance circle packings on surfaces}
\author{Xu Xu, Chao Zheng}

\date{\today}

\address{School of Mathematics and Statistics, Wuhan University, Wuhan, 430072, P.R.China}
 \email{xuxu2@whu.edu.cn}

\address{School of Mathematics and Statistics, Wuhan University, Wuhan 430072, P.R. China}
\email{czheng@whu.edu.cn}

\thanks{MSC (2020): 52C25, 52C26.}

\keywords{Combinatorial curvature flows; Surgery; Inversive distance circle packings; Piecewise Euclidean metrics}

\begin{abstract}
Inversive distance circle packings introduced by Bowers-Stephenson  are
natural generalizations of Thurston's circle packings on surfaces.
To find piecewise Euclidean metrics on surfaces with prescribed combinatorial curvatures,
we introduce the combinatorial Calabi flow, the fractional combinatorial Calabi flow and the combinatorial $p$-th Calabi flow
for the Euclidean inversive distance circle packings.
Due to the singularities possibly developed by these combinatorial curvature flows,
the longtime existence and convergence of these combinatorial curvature flows have been a difficult problem for a long time.
To handle the potential singularities along these combinatorial curvature flows, we do surgery along these flows by edge flipping
under the weighted Delaunay condition.
Using the discrete conformal theory recently established by Bobenko-Lutz for decorated piecewise Euclidean metrics on surfaces,
we prove the longtime existence and global convergence for the solutions of these combinatorial curvature flows with surgery.
This provides effective algorithms for finding piecewise Euclidean metrics on surfaces with prescribed combinatorial curvatures.
\end{abstract}

\maketitle

\section{Introduction}
\subsection{PE metrics and PE surfaces}

Suppose $S$ is a connected closed surface  and $V$ is a finite subset of $S$ with $|V|=N$, we call $(S, V)$  a marked surface.
A piecewise Euclidean (PE for short) metric $dist_{S}$ on the marked surface $(S,V)$ is a flat cone metric with the conic singularities contained in $V$.
The marked surface $(S, V)$ endowed with a PE metric $dist_{S}$ is called a PE surface.
Let $\mathcal{T}={(V,E,F)}$ be a triangulation of $(S, V)$, where $V,E,F$ are the sets of vertices, edges and faces respectively.
We use one index to denote a vertex (such as $i$), two indices to denote an edge (such as $\{ij\}$) and three indices to denote a face (such as $\{ijk\}$)
in the triangulation $\mathcal{T}$.
The marked surface $(S, V)$ with a fixed triangulation $\mathcal{T}$ is called a triangulated surface, denoted by $(S, \mathcal{T})$.
The triangulation $\mathcal{T}$ for a PE surface is a geodesic triangulation if the edges are geodesics in the PE metric.
The PE metric $dist_{S}$ on a geodesic triangulated PE surface $(S,\mathcal{T}, dist_{S})$ defines a map $l: E\rightarrow \mathbb{R}_{>0}$ such that $l_{ij}, l_{ik}, l_{jk}$ satisfy the triangle inequalities for any triangle $\{ijk\}\in F$.
Conversely, given a function $l: E\rightarrow \mathbb{R}_{>0}$ satisfying the triangle inequality,
one can construct a PE metric on a triangulated surface $(S,\mathcal{T})$ by isometrically gluing Euclidean triangles along edges in pairs.
Therefore, we use $l: E\rightarrow \mathbb{R}_{>0}$ to denote a PE metric on a triangulated surface $(S,\mathcal{T})$.
For a PE metric $l: E\rightarrow \mathbb{R}_{>0}$ on a triangulated surface $(S,\mathcal{T})$,
the combinatorial curvature $K: V\rightarrow (-\infty, 2\pi)$ is used to describe the conic singularities of PE metrics at the vertices.
For a vertex $i\in V$, the combinatorial curvature $K_i$ is defined to be
\begin{equation*}
K_i=2\pi-\sum_{\{ijk\}\in F}\theta^i_{jk},
\end{equation*}
where the summation is taken over all the triangles with $i$ as a vertex and $\theta^i_{jk}$ is the inner angle of the triangle $\{ijk\}$ at the vertex $i$.
Note that the combinatorial curvature $K$ is independent of the geodesic triangulations of a PE surface.
Therefore, the combinatorial curvature $K$ is an intrinsic geometric invariant of a PE surface.
Furthermore, the combinatorial curvature $K$ satisfies the following discrete Gauss-Bonnet formula (\cite{Chow-Luo}, Proposition 3.1)
\begin{equation}\label{Eq: Gauss-Bonnet}
\sum_{i=1}^N K_i=2\pi \chi(S),
\end{equation}
where $\chi(S)$ is the Euler number of the surface $S$.

Due to the fast developing of computer science, applied mathematics, engineering and medical imaging,
PE surfaces with prescribed combinatorial curvatures have been more and more important in theory and applications.
Please refer to \cite{ZG,Gu-Yau} and others for this.
A basic problem in theory and applications is finding PE metrics on surfaces with prescribed combinatorial curvatures.
An effective approach to this problem is studying it in discrete conformal geometry and finding the PE metrics with prescribed combinatorial curvatures by combinatorial curvature flows.
There are mainly two types of discrete conformal structures on surfaces that have been extensively
studied in discrete conformal geometry in the literature.
One type is the vertex scalings introduced by Luo \cite{L1}.
Luo's vertex scaling of a PE metric $l$ on a triangulated surface $(S,\mathcal{T})$ is defined to be the PE metric $\widetilde{l}$ on $(S,\mathcal{T})$ such that
there exists a discrete conformal factor $u\in \mathbb{R}^V$ with
\begin{equation}\label{Eq: VS}
\widetilde{l}_{ij}=l_{ij}e^{\frac{u_i+u_j}{2}}, \ \ \forall \{ij\}\in E.
\end{equation}
In the last two decades, there have been lots of important works on Luo's vertex scalings.
The rigidity of Luo's vertex scalings was proved locally by Luo \cite{L1} and globally by Bobenko-Pinkall-Springborn \cite{BPS}.
One can also refer to \cite{XZ1} for an elementary proof of the rigidity of Luo's vertex scalings.
The discrete uniformization theorems for Luo's vertex scalings were established by Gu-Luo-Sun-Wu \cite{Gu1}, Gu-Guo-Luo-Sun-Wu \cite{Gu2}, Springborn \cite{Springborn} and Izmestiev-Prosanov-Wu \cite{IPW}.
The convergence of Luo's vertex scalings was studied by Gu-Luo-Wu \cite{GLW}, Luo-Sun-Wu \cite{LSW} and Wu-Zhu \cite{WZ}.
To find PE metrics with prescribed combinatorial curvatures on surfaces,
the combinatorial curvature flows for Luo's vertex scalings were extensively studied.
Luo \cite{L1} first introduced the combinatorial Yamabe flow for the vertex scalings on a triangulated surface and proved its local convergence.
Ge-Jiang \cite{GJ0} further proved the global convergence of the combinatorial Yamabe flow for Luo's vertex scalings on a triangulated surface
by the constant extension introduced by Bobenko-Pinkall-Springborn \cite{BPS}.
The combinatorial Yamabe flows with surgery for Luo's vertex scalings were introduced by Gu-Luo-Sun-Wu \cite{Gu1} and Gu-Guo-Luo-Sun-Wu \cite{Gu2},
where the longtime existence and convergence were also proved.
The finiteness of surgeries along the combinatorial Yamabe flow with surgery was proved by Wu \cite{Wu}.
Following Luo's work on combinatorial Yamabe flow \cite{L1},
Ge \cite{Ge1} introduced the combinatorial Calabi flow for Luo's vertex scalings on triangulated surfaces and proved the short time existence.
The longtime existence and convergence of the combinatorial Calabi flow with surgery for Luo's vertex scalings were proved by Zhu-Xu \cite{ZX}.
Feng-Lin-Zhang \cite{FLZ} introduced the combinatorial $p$-th Calabi flow for Luo's vertex scalings on surfaces,
and proved the corresponding longtime existence and convergence of the combinatorial $p$-th Calabi flow by surgery.
Recently, Wu-Xu \cite{Wu-Xu} introduced the fractional combinatorial Calabi flow for Luo's vertex scalings on surfaces, and proved the corresponding longtime existence and convergence of the flow by surgery.
This generalizes the results for the combinatorial Yamabe flow and the combinatorial Calabi flow with surgery in
\cite{Gu1, Gu2, ZX}.
Another important type of discrete conformal structures on surfaces is the circle packings.
Among different types of circle packings, Thurston's circle packings have been extensively studied in the literature, including the existence, rigidity, convergence and the combinatorial curvature flows.
However, not all PE metrics are Thurston's circle packing metrics.
The inversive distance circle packings introduced by Bowers-Stephenson \cite{BS} are natural generalizations of Thurston's circle packings.
Comparing to Luo's vertex scalings, the corresponding theory for the inversive distance circle packings is not well established.

\subsection{Inversive distance circle packings and decorated PE metrics}
The inversive distance circle packing metrics on surfaces are defined as follows.
\begin{definition}[\cite{BS}]\label{Def: IDCP}
Suppose $(S,\mathcal{T})$ is a triangulated surface with a weight $I: E\rightarrow (-1, +\infty)$.
A PE metric $l: E\rightarrow \mathbb{R}_{>0}$ on $(S,\mathcal{T},I)$ is an inversive distance circle packing metric
if there exists a function $r: V \rightarrow\mathbb{R}_{>0}$ such that
\begin{equation}\label{Eq: IDCP}
l_{ij}=\sqrt{r^2_i+r^2_j+2 I_{ij}r_ir_j}
\end{equation}
for any edge $\{ij\}\in E$.
\end{definition}
The map $r: V \rightarrow\mathbb{R}_{>0}$ is referred as an \textit{inversive distance circle packing} on $(S, \mathcal{T}, I)$.
The map $r$ is also called a decoration on $(S,\mathcal{T},l)$ by Bobenko-Lutz \cite{BL}.
The pair $(l,r)$ on a triangulated surface $(S,\mathcal{T})$ is called a decorated PE metric.
The weight $I_{ij}$ is the inversive distance of
the two circles centered at $i$ and $j$ with radii $r_i$ and $r_j$ respectively.
If $I_{ij}\in(-1,0)$, then the two circles attached to the vertices $i$ and $j$ intersect with an obtuse angle.
If $I_{ij}\in[0,1]$, then the two circles intersect with a non-obtuse angle.
Taking $I_{ij}=\cos \Phi_{ij}$ with $\Phi_{ij}\in [0, \frac{\pi}{2}]$, the inversive distance circle packings are reduced to Thurston's circle packings in \cite{Thurston}.
If $I_{ij}\in(1,+\infty)$, then the two circles are disjoint.

Note that the inversive distance between two circles is invariant under M\"{o}bius transformations \cite{Coxeter}.
Two inversive distance circle packing metrics
$l$ and $\widetilde{l}$ on $(S, \mathcal{T}, I)$ and $(S, \mathcal{T}, \widetilde{I})$ respectively
are discrete conformally equivalent if $I=\widetilde{I}$.
In this case, we set
\begin{equation}\label{Eq: DCE1}
\widetilde{r}_i=e^{u_i}r_i
\end{equation}
for $i\in V$ and call $u: V\rightarrow \mathbb{R}$ a discrete conformal factor.
Then the equation (\ref{Eq: IDCP}) is equivalent to
\begin{equation}\label{Eq: DCE2}
\widetilde{l}_{ij}^2
=(e^{2u_i}-e^{u_i+u_j})r^2_i
+(e^{2u_j}-e^{u_i+u_j})r^2_j
+e^{u_i+u_j}l_{ij}^2
\end{equation}
for any edge $\{ij\}\in E$.
Conversely, if two decorated PE metrics $(l,r)$ and $(\widetilde{l},\widetilde{r})$ on $(S,\mathcal{T})$ satisfy
the equations (\ref{Eq: DCE1}) and (\ref{Eq: DCE2}) for any edge $\{ij\}\in E$, then they are discrete conformallly equivalent.
Please refer to Bobenko-Lutz's recent work \cite{BL} for more information on this.

\begin{remark}\label{Rmk: 1}
By comparing (\ref{Eq: VS}) with (\ref{Eq: DCE2}),
one can see that they can be written in a unified form
\begin{equation}\label{Eq: DCS}
\widetilde{l}^2_{ij}=\alpha_ie^{2f_i}+\alpha_je^{2f_j}
+2\eta_{ij}e^{f_i+f_j}
\end{equation}
with $\alpha_i: V\rightarrow \mathbb{R}$ and $\eta_{ij}: E\rightarrow \mathbb{R}$.
The equation (\ref{Eq: DCS}) characterizes the discrete conformal structures introduced by
Glickenstein \cite{Glickenstein2011}, Glickenstein-Thomas \cite{GT} and Zhang-Guo-Zeng-Luo-Yau-Gu \cite{ZGZLYG}.
If we set $\alpha_i=0$ and $2\eta_{ij}=l_{ij}^2$, then (\ref{Eq: DCS}) is equivalent to (\ref{Eq: VS}).
If we set $\alpha_i=r_i^2>0$ and $2\eta_{ij}=l_{ij}^2-r^2_i-r^2_j$, then (\ref{Eq: DCS}) is equivalent to (\ref{Eq: DCE2}).
\end{remark}

Bowers-Stephenson \cite{BS} conjectured that the inversive distance circle packings on surfaces are rigid.
This conjecture was proved to be true by Guo \cite{Guo}, Luo \cite{L3} and Xu \cite{Xu1, Xu2} for the inversive distance in $(-1,+\infty)$.
For the inversive distance in $(1,+\infty)$,
the convergence of inversive distance circle packings was studied by Chen-Luo-Xu-Zhang \cite{CLXZ}.
Recently, Bobenko-Lutz \cite{BL} introduced a new definition of discrete conformality for Euclidean inversive distance circle packings, i.e., decorated PE metrics, and proved the corresponding discrete uniformization theorem.
However, there are not so much research on the combinatorial curvature flows for inversive distance circle packings on surfaces.

\subsection{Combinatorial curvature flows for inversive distance circle packings and the main results}
To find effective algorithms searching for polyhedral metrics with prescribed combinatorial curvatures on surfaces,
Chow-Luo \cite{Chow-Luo} introduced the combinatorial Ricci flow.
Motivated by Chow-Luo's original work \cite{Chow-Luo},
Ge \cite{Ge1,Ge2} and Ge-Xu \cite{GX DGA} introduced the combinatorial Calabi flow,
Lin-Zhang \cite{L-Z} introduced the combinatorial $p$-th Calabi flow
and Wu-Xu \cite{Wu-Xu} introduced the fractional combinatorial Calabi flow on surfaces.
These combinatorial curvature flows were initially introduced for Thurston's circle packings.
The longtime existence and convergence of these combinatoiral curvature flows
for Thurston's circle packings were proved in \cite{Chow-Luo, Ge1,Ge2,GH, GX DGA,Wu-Xu,L-Z}.
However, not all PE metrics are Thurston's circle packing metrics.
As natural generalizations of Thurston's circle packings,
the inversive distance circle packings have more flexibility than Thurston's circle packings.
On the other hand, due to the possible degenerations of triangles generated by inversive distance circle packings,
the longtime existence and convergence of the solutions of the combinatorial curvature flows for inversive distance circle packings
have been a difficult problem for a long time.
The combinatorial Ricci flow for inversive distance circle packings on triangulated surfaces was first introduced by Zhang-Guo-Zeng-Luo-Yau-Gu \cite{ZGZLYG}.
The properties of the combinatorial Ricci flow for inversive distance circle packings on triangulated surfaces were studied in \cite{GJ1, GJ2, GJ3, GX JFA} by the constant extension introduced by Luo \cite{L3}.
Recently, Bobenko-Lutz \cite{BL} proved the longtime existence and convergence of the combinatorial Ricci flow with surgery
for the Euclidean inversive distance circle packings on surfaces with inversive distance in $(1,+\infty)$.
The efficiency, efficacy and robustness of this combinatorial Ricci flow with surgery in applications were demonstrated in \cite{CZLG}.
In this paper, we introduce the combinatorial Calabi flow, the fractional combinatorial Calabi flow and the combinatorial $p$-th Calabi flow for
the Euclidean inversive distance circle packings on surfaces with inversive distance in $(1,+\infty)$.
Furthermore, we prove the longtime existence and convergence of the solutions of these combinatorial curvature flows by surgery.

\begin{definition}
Suppose $(S,\mathcal{T})$ is a triangulated surface with a weight $I: E\rightarrow (1, +\infty)$.
Let $\overline{K}: V\rightarrow (-\infty, 2\pi)$ be a given function defined on $V$.
The combinatorial Calabi flow for the inversive distance circle packings on $(S,\mathcal{T}, I)$ is defined to be
\begin{eqnarray}\label{Eq: CCF}
\begin{cases}
\frac{du_i}{dt}=\Delta_{\mathcal{T}}(K-\overline{K})_i,\\
u_i(0)=u_0,
\end{cases}
\end{eqnarray}
where $\Delta_{\mathcal{T}}$ is the discrete Laplace operator defined by
\begin{equation}\label{Eq: Laplace}
\Delta_{\mathcal{T}}f_i=-\sum_{j=1}^N\frac{\partial K_i}{\partial u_j}f_j
\end{equation}
for any function $f: V\rightarrow \mathbb{R}$ defined on the vertices.
\end{definition}

\begin{remark}
The combinatorial Calabi flow (\ref{Eq: CCF})
is a negative gradient flow of the combinatoiral Calabi energy
$$\mathcal{C}(u):=||K-\overline{K}||^2=\sum_{i=1}^N (K_i(u)-\overline{K}_i)^2.$$
This fact was first observed by Ge \cite{Ge1,Ge2} for Thurston's Euclidean circle packings on surfaces.
\end{remark}

Set
\begin{equation*}
L=(L_{ij})_{N\times N}
=\frac{\partial(K_1,...,K_N)}{\partial(u_1,...,u_N)}.
\end{equation*}
The equation (\ref{Eq: Laplace}) implies $\Delta_{\mathcal{T}}=-L$.
By Lemma \ref{Lem: Xu}, the matrix $L$ is symmetric and positive semi-definite on the admissible space $\Omega^{\mathcal{T}}$ of discrete conformal factors $u$ such that the equation (\ref{Eq: DCE2}) defines a PE metric.
Then there exists an orthonormal matrix $P$ such that
\begin{equation*}
L=P^\mathrm{T}\cdot \text{diag}\{\lambda_1,...,\lambda_N\}\cdot P,
\end{equation*}
where $\lambda_1,...,\lambda_N$ are non-negative eigenvalues of the matrix $L$.
For any $s\in \mathbb{R}$, the $2s$-th order fractional discrete Laplace operator $\Delta_{\mathcal{T}}^s$ is defined to be
\begin{equation}\label{Eq: fractional Laplace}
\Delta_{\mathcal{T}}^s=-L^s=-P^\mathrm{T}\cdot \text{diag}\{\lambda^s_1,...,\lambda^s_N\}\cdot P.
\end{equation}
Therefore, the $2s$-th order fractional discrete Laplace operator $\Delta_{\mathcal{T}}^s$ is negative semi-definite on $\Omega^{\mathcal{T}}$.
Furthermore, the $2s$-th order fractional discrete Laplace operator $\Delta_{\mathcal{T}}^s$
has the same kernel space as the discrete Laplace operator $\Delta_{\mathcal{T}}$.
Specially, if $s=0$, then $\Delta_{\mathcal{T}}^s$ is reduced to the minus identity operator;
if $s=1$, then $\Delta_{\mathcal{T}}^s$ is reduced to the discrete Laplace operator $\Delta_{\mathcal{T}}=-L=-(\frac{\partial K_i}{\partial u_j})_{N\times N}$.

\begin{definition}
Suppose $(S,\mathcal{T})$ is a triangulated surface with a weight $I: E\rightarrow (1, +\infty)$.
Let $\overline{K}: V\rightarrow (-\infty, 2\pi)$ be a given function defined on $V$.
For any $s\in \mathbb{R}$,
the $2s$-th order fractional combinatorial Calabi flow for the inversive distance circle packings on $(S,\mathcal{T}, I)$ is defined to be
\begin{eqnarray}\label{Eq: FCCF}
\begin{cases}
\frac{du_i}{dt}=\Delta_{\mathcal{T}}^s(K-\overline{K})_i,\\
u_i(0)=u_0,
\end{cases}
\end{eqnarray}
where $\Delta_{\mathcal{T}}^s$ is the fractional discrete Laplace operator defined by (\ref{Eq: fractional Laplace}).
\end{definition}

\begin{remark}
If $s=0$, the $2s$-th order fractional combinatorial Calabi flow (\ref{Eq: FCCF}) is reduced to the combinatorial Ricci flow introduced by Zhang-Guo-Zeng-Luo-Yau-Gu \cite{ZGZLYG} for the inversive distance circle packings.
If $s=1$, the $2s$-th order fractional combinatorial Calabi flow (\ref{Eq: FCCF}) is reduced to the combinatorial Calabi flow (\ref{Eq: CCF}).
\end{remark}

By Lemma \ref{Lem: Xu}, we have $\sum_{j=1}^N\frac{\partial K_i}{\partial u_j}=0$.
As a result, the equation (\ref{Eq: Laplace}) defining the discrete Laplace operator $\Delta_{\mathcal{T}}$ can be written as
\begin{equation*}
\Delta_{\mathcal{T}}f_i
=-\sum_{j=1}^N\frac{\partial K_i}{\partial u_j}f_j
=\sum_{j\sim i}(-\frac{\partial K_i}{\partial u_j})(f_j-f_i).
\end{equation*}
For any $p>1$, we define the discrete $p$-th Laplace operator $\Delta_{p,\mathcal{T}}$ for the inversive distance circle packings on
$(S,\mathcal{T}, I)$ by the following formula
\begin{equation}\label{Eq: P-Laplace}
\Delta_{p,\mathcal{T}}f_i=\sum_{j\sim i}(-\frac{\partial K_i}{\partial u_j})|f_j-f_i|^{p-2}(f_j-f_i),
\end{equation}
where $f: V\rightarrow \mathbb{R}$ is a function defined on the vertices.

\begin{definition}\label{Def: PCCF}
Suppose $(S,\mathcal{T})$ is a triangulated surface with a weight $I: E\rightarrow (1, +\infty)$.
Let $\overline{K}: V\rightarrow (-\infty, 2\pi)$ be a given function defined on $V$.
For any $p>1$,
the combinatorial $p$-th Calabi flow for the inversive distance circle packings on $(S,\mathcal{T}, I)$ is defined to be
\begin{eqnarray}\label{Eq: PCCF}
\begin{cases}
\frac{du_i}{dt}=\Delta_{p,\mathcal{T}}(K-\overline{K})_i,\\
u_i(0)=u_0,
\end{cases}
\end{eqnarray}
where $\Delta_{p,\mathcal{T}}$ is the discrete $p$-th Laplace operator defined by (\ref{Eq: P-Laplace}).
\end{definition}

\begin{remark}
If $p=2$, the discrete $p$-th Laplace operator $\Delta_{p,\mathcal{T}}$ is reduced to the discrete Laplace operator $\Delta_{\mathcal{T}}$ and the combinatorial $p$-th Calabi flow (\ref{Eq: PCCF}) is reduced to the combinatorial Calabi flow (\ref{Eq: CCF}).
\end{remark}

As the combinatorial Calabi flow (\ref{Eq: CCF}), the $2s$-th order fractional combinatorial Calabi flow (\ref{Eq: FCCF}) and the combinatorial $p$-th Calabi flow (\ref{Eq: PCCF}) are ODE systems with smooth coefficients,
the solutions of these combinatorial curvature flows always exist locally around the initial time $t=0$.
Furthermore, we have the following results on the longtime existence and convergence for the solutions of the combinatorial Calabi flow (\ref{Eq: CCF}) and the $2s$-th order fractional combinatorial Calabi flow (\ref{Eq: FCCF}).

\begin{theorem}\label{Thm: main 1}
Suppose $(S,\mathcal{T})$ is a triangulated surface with a weight $I: E\rightarrow (1, +\infty)$.
Let $\overline{K}: V\rightarrow (-\infty, 2\pi)$ be a given function defined on $V$ satisfying the discrete Gauss-Bonnet formula (\ref{Eq: Gauss-Bonnet})
and $s\in \mathbb{R}$ be a constant.
\begin{description}
\item[(i)] If the solution $u(t)$ of the combinatorial Calabi flow (\ref{Eq: CCF}) or the $2s$-th order fractional combinatorial Calabi flow (\ref{Eq: FCCF}) converges, then there exists a discrete conformal factor on $(S,\mathcal{T}, I)$ with the combinatorial curvature $\overline{K}$.
  \item[(ii)] If there exists a discrete conformal factor $\overline{u}$ with the combinatorial curvature $\overline{K}$ on $(S,\mathcal{T}, I)$, there exists a constant $\delta>0$ such that if the initial value $u(0)$ satisfies $||u(0)-\overline{u}||<\delta$ and $\sum^N_{i=1}u_i(0)=\sum^N_{i=1}\overline{u}_i$,
then the solutions of the combinatorial Calabi flow (\ref{Eq: CCF}) and the $2s$-th order fractional combinatorial Calabi flow (\ref{Eq: FCCF}) exist for all time and converge exponentially fast to $\overline{u}$.
\end{description}
\end{theorem}

\begin{figure}[!ht]
\centering
\begin{overpic}[scale=0.3]{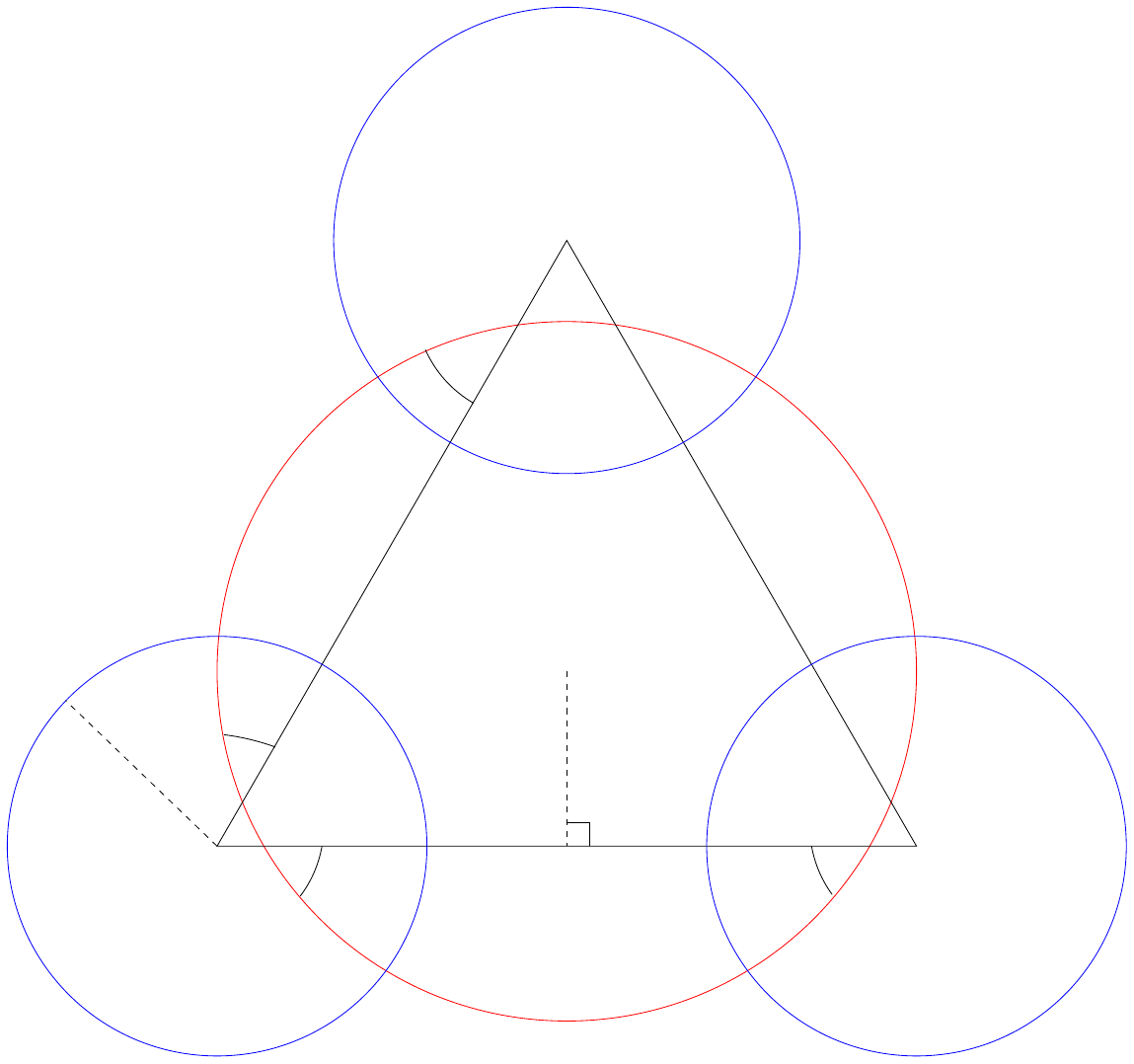}
\put (15,15) {$i$}
\put (80.5,15) {$j$}
\put (48.5,74) {$k$}
\put (30,55) {$\alpha_{ik}^j$}
\put (20,30) {$\alpha_{ik}^j$}
\put (27,12.5) {$\alpha_{ij}^k$}
\put (63,12.5) {$\alpha_{ij}^k$}
\put (10,25) {$r_i$}
\put (51,25) {$d_{ij}^k$}
\put (48,13) {$r_{ij}$}
\put (48,38) {$c_{ijk}$}
\end{overpic}
\caption{Datas for a decorated triangle $\{ijk\}\in F$}
\label{figure}
\end{figure}

For general initial inversive distance circle packing metrics, the combinatorial Calabi flow (\ref{Eq: CCF}), the $2s$-th order fractional combinatorial Calabi flow (\ref{Eq: FCCF}) and the combinatorial $p$-th Calabi flow (\ref{Eq: PCCF}) may develop singularities, including that some triangles degenerate and the discrete conformal factors tend to infinity along these combinatorial curvature flows.
Specially, it is proved (\cite{Xu2}, Remark 2.6) that $\frac{\partial \theta^i_{jk}}{\partial u_j}$ tends to infinity if the nondegenerate triangle generated by inversive distance circle packings tends to be a degenerate triangle (the triangle inequality fails).
To handle the potential singularities along these combinatorial curvature flows,
we do surgery on the flows by edge flipping under the weighted Delaunay condition, the idea of which comes from the recent work of Bobenko-Lutz \cite{BL}.
Given a decorated triangle $\{ijk\}\in F$, there is a unique circle $C_{ijk}$ simultaneously orthogonal to all the three circles attached to the vertices $i,j,k$.
The circle  $C_{ijk}$ is called the face-circle of the decorated triangle $\{ijk\}$.
Denote $\alpha_{ij}^k$ as the interior intersection angle of the face-circle $C_{ijk}$ and the edge $\{ij\}$.
Please refer to Figure \ref{figure}.
A \emph{weighted Delaunay triangulation} $\mathcal{T}$ for a decorated PE metric $(l,r)$ on $(S,V)$ is a geodesic triangulation such that $\alpha_{ij}^k+\alpha_{ij}^l\leq \pi$ for any adjacent triangles $\{ijk\}$ and $\{ijl\}$ sharing a common edge $\{ij\}\in E$.
This definition of weighted Delaunay triangulations comes from Bobenko-Lutz \cite{BL}.
One can also refer to \cite{Glickenstein DCG, Glickenstein arxiv, Xu 21} for other equivalent definitions of weighted Delaunay triangulations.
Along these combinatorial curvature flows (\ref{Eq: CCF}),  (\ref{Eq: FCCF}) and (\ref{Eq: PCCF})
on $(S,\mathcal{T})$, if $\mathcal{T}$ is weighted Delaunay in $(l(t),r(t))$ for $t\in [0,T]$ and not weighted Delaunay in $(l(t),r(t))$ for $t\in (T,T+\epsilon),\ \epsilon>0$, there exists an edge $\{ij\}\in E$ such that $\alpha_{ij}^k+\alpha_{ij}^l\leq \pi$ for $t\in [0,T]$ and $\alpha_{ij}^k+\alpha_{ij}^l> \pi$ for $t\in (T,T+\epsilon)$.
Then we replace the triangulation $\mathcal{T}$ by a new triangulation $\mathcal{T}'$ at time $t=T$ via replacing two triangles $\{ijk\}$ and $\{ijl\}$ adjacent to $\{ij\}$ by two new triangles $\{ikl\}$ and $\{jkl\}$. This procedure is called a \textbf{surgery\  by\  flipping} on the triangulation $\mathcal{T}$, which is also an isometry of $(S,V)$ with the decorated PE metric $(l(T),r(T))$.
After the surgery by flipping at the time $t=T$, we run these combinatorial curvature flows on $(S, \mathcal{T}')$
with the initial metric $(l(T),r(T))$.
Whenever the weighted Delaunay condition is not satisfied along these combinatorial curvature flows, we do surgery on these combinatorial curvature flows by flipping.

We have the following result on the longtime existence and convergence for the solutions of the combinatorial Calabi flow with surgery, the $2s$-th order fractional combinatorial Calabi flow with surgery and the combinatorial $p$-th Calabi flow with surgery for decorated PE metrics on $(S,V)$.

\begin{theorem}\label{Thm: main 2}
Suppose $(S,V)$ is a marked surface with a decorated PE metric $(dist_g,r)$.
Let $\overline{K}: V\rightarrow (-\infty, 2\pi)$ be a given function defined on $V$ satisfying the discrete Gauss-Bonnet formula (\ref{Eq: Gauss-Bonnet}).
\begin{description}
\item[(i)] The solution of the combinatorial Calabi flow with surgery exists for all time and converges exponentially fast to $\overline{u}$ with the prescribed combinatorial curvature $\overline{K}$ for any initial $u(0)\in \mathbb{R}^N$ with $\sum^N_{i=1}u_i(0)=\sum^N_{i=1}\overline{u}_i$;
\item[(ii)] For any $s\in \mathbb{R}$, the solution of the $2s$-th order fractional combinatorial Calabi flow with surgery exists for all time and converges exponentially fast to $\overline{u}$ with the prescribed combinatorial curvature $\overline{K}$ for any initial $u(0)\in \mathbb{R}^N$ with $\sum^N_{i=1}u_i(0)=\sum^N_{i=1}\overline{u}_i$;
\item[(iii)] For any $p>1$, the solution of the combinatorial $p$-th Calabi flow with surgery exists for all time and converges to $\overline{u}$ with the prescribed combinatorial curvature $\overline{K}$ for any initial $u(0)\in \mathbb{R}^N$ with $\sum^N_{i=1}u_i(0)=\sum^N_{i=1}\overline{u}_i$.
\end{description}
\end{theorem}

\begin{remark}
If $s=0$, the convergence for the $2s$-th order fractional combinatorial Calabi flow with surgery in Theorem \ref{Thm: main 2} is reduced to the convergence
of combinatorial Ricci flow with surgery for decorated PE metrics obtained by Bobenko-Lutz \cite{BL}.
Different from the combinatorial Calabi flow with surgery ($p=2$),
we can not get the exponential convergence for the solution of the combinatorial $p$-th Calabi flow with surgery for $p\neq2$.
\end{remark}

In applications, it is important that we just need to do finite times of surgeries along these combinatorial curvature flows.
This is closely related to the stability of computer algorithms.
In the case of Luo's vertex scalings, Gu-Luo-Sun-Wu \cite{Gu1} proved that $\mathbb{R}^N=\bigcup_{i=1}^m\mathcal{D}_i$ is an analytical decomposition of cells.
Hence, this problem is reduced to whether the solutions of the combinatorial curvature flows with surgery go cross the boundaries of these cells only finitely many times.
Wu \cite{Wu} solved this problem for the combinatorial Yamabe flow of Luo's vertex scalings and proved the following theorem.
\begin{theorem}[\cite{Wu}]
Suppose $\mathbb{R}^N=\bigcup_{i=1}^m\mathcal{D}_i$ is an analytic cell decomposition.
Let $f(x)\in C^1(\mathbb{R}^N)$ be analytic on each cell $\mathcal{D}_i$ and have a unique minimum point where $f$ has positive Hession.
Then its gradient flow $\gamma(t)$ satisfying $\gamma^\prime(t)=-\nabla f(\gamma(t))$ intersects the cell faces $\mathcal{D}_i$
finitely many times.
\end{theorem}

For decorated PE metrics on surfaces,
Bobenko-Lutz \cite{BL} recently proved similar analytic cell decomposition of $\mathbb{R}^N$.
It is convinced that we also just need to do finite times of surgeries along the combinatorial Calabi flow, the $2s$-th order fractional combinatorial Calabi flow and the combinatorial $p$-th Calabi flow for decorated PE metrics on $(S,V)$.

\subsection{Organization of the paper}
The paper is organized as follows.
In Section \ref{Sec: fixed triangulation}, we study the combinatorial Calabi flow (\ref{Eq: CCF}), the $2s$-th order fractional combinatorial Calabi flow (\ref{Eq: FCCF}) and the combinatorial $p$-th Calabi flow (\ref{Eq: PCCF}) on triangulated surfaces and prove Theorem \ref{Thm: main 1}.
In Section \ref{Sec: variable triangulation}, we allow the triangulation on a marked surface to be changed by edge flipping and prove Theorem \ref{Thm: main 2}.
\\
\\
\textbf{Acknowledgements}\\[8pt]
The first author thanks Professor Feng Luo for his invitation to the workshop
``Discrete and Computational Geometry, Shape Analysis, and Applications" taking place
at Rutgers University, New Brunswick from May 19th to May 21st, 2023.
The first author also thanks Carl O. R. Lutz for helpful communications during the workshop.

\section{Combinatorial curvature flows for fixed triangulations}\label{Sec: fixed triangulation}
Suppose $(S,\mathcal{T},I)$ is a weighted triangulated surface with a decorated PE metric $(l,r)$.
The admissible space $\Omega_{ijk}^{\mathcal{T}}$ of the discrete conformal factors for a triangle $\{ijk\}\in F$ in $(S,\mathcal{T})$ is defined to be the set of $({u_i,u_j,u_k})\in \mathbb{R}^3$ such that the triangle with edge lengths
$\widetilde{l}_{ij}, \widetilde{l}_{ik}, \widetilde{l}_{jk}$ defined by (\ref{Eq: DCE2})
exists in the 2-dimensional Euclidean space $\mathbb{E}^2$, i.e.,
\begin{equation*}
\Omega^{\mathcal{T}}_{ijk}=\{(u_i,u_j,u_k)\in \mathbb{R}^3| \widetilde{l}_{rs}+\widetilde{l}_{rt}>\widetilde{l}_{ts}, \{r,s,t\}=\{i,j,k\} \}.
\end{equation*}
The admissible space of discrete conformal factors on $(S,\mathcal{T},I)$, denoted by $\Omega^{\mathcal{T}}$,
is defined to be the vectors $u\in \mathbb{R}^N$ such that  $({u_i,u_j,u_k})\in \Omega^{\mathcal{T}}_{ijk}$ for every triangle $\{ijk\}\in F$.

\begin{lemma}[\cite{Guo, Xu1,Xu2}]\label{Lem: Xu}
Suppose $(S,\mathcal{T})$ is a triangulated surface with a weight $I: E\rightarrow (1, +\infty)$.
\begin{description}
  \item[(i)] The admissible space $\Omega^{\mathcal{T}}_{ijk}$ is a non-empty simply connected open set whose boundary is analytic.
  \item[(ii)] The matrix $\frac{\partial (\theta^{i}_{jk}, \theta^{j}_{ik}, \theta^{k}_{ij})}{\partial (u_i, u_j, u_k)}$ is symmetric and negative semi-definite with rank 2 and kernel $\{c(1,1,1)^\mathrm{T}|c\in \mathbb{R} \}$ on $\Omega^{\mathcal{T}}_{ijk}$.
       As a result, the matrix $L=\frac{\partial(K_1,...,K_N)}{\partial(u_1,...,u_N)}$ is symmetric and positive semi-definite with rank $N-1$ and kernel $\{c\mathbf{1}^\mathrm{T}\in \mathbb{R}^N|c\in \mathbb{R}\}$ on $\Omega^\mathcal{T}$.
\end{description}
\end{lemma}

By Lemma \ref{Lem: Xu}, we have the following result.
\begin{lemma}\label{Lem: invariant 1}
Suppose $(S,\mathcal{T})$ is a triangulated surface with a weight $I: E\rightarrow (1, +\infty)$.
Let $\overline{K}: V\rightarrow (-\infty, 2\pi)$ be a given function defined on $V$.
If $\overline{K}$ satisfies the discrete Gauss-Bonnet formula (\ref{Eq: Gauss-Bonnet}),
then $\sum_{i=1}^Nu_i(t)$ is invariant along the combinatorial Calabi flow (\ref{Eq: CCF}), the $2s$-th order fractional combinatorial Calabi flow (\ref{Eq: FCCF}) and the combinatorial $p$-th Calabi flow (\ref{Eq: PCCF}).
\end{lemma}
\proof
By Lemma \ref{Lem: Xu} and direct calculations, we have
\begin{equation*}
\frac{d(\sum_{i=1}^N u_i)}{dt}
=\sum_{i=1}^N\Delta_{\mathcal{T}}(K-\overline{K})_i
=-\sum_{i=1}^N\sum_{j=1}^N(\frac{\partial K_i}{\partial u_j})(K-\overline{K})_j
=0
\end{equation*}
along the combinatorial Calabi flow (\ref{Eq: CCF}).
This implies $\sum_{i=1}^N u_i$ is invariant along the combinatorial Calabi flow (\ref{Eq: CCF}).
Similarly, by Lemma \ref{Lem: Xu}, we have
\begin{equation*}
\frac{d(\sum_{i=1}^N u_i)}{dt}
=\sum_{i=1}^N\Delta^s_{\mathcal{T}}(K-\overline{K})_i
=-\mathbf{1}^\mathrm{T}(\frac{\partial K}{\partial u})^s(K-\overline{K})
=0
\end{equation*}
along the $2s$-th order fractional combinatorial Calabi flow (\ref{Eq: FCCF}).
This implies $\sum_{i=1}^N u_i$ is invariant along the $2s$-th order fractional combinatorial Calabi flow (\ref{Eq: FCCF}).
Similarly, by direct calculations, we have
\begin{equation*}
\frac{d(\sum_{i=1}^N u_i)}{dt}
=\sum_{i=1}^N\Delta_{p,\mathcal{T}}(K-\overline{K})_i
\end{equation*}
along the combinatorial $p$-th Calabi flow (\ref{Eq: PCCF}).
By the following formula obtained by Lin-Zhang (\cite{L-Z}, Lemma 5.3)
\begin{equation*}
\sum_{i=1}^N\Delta_{p,\mathcal{T}}f_i=0
\end{equation*}
for any $f: V\rightarrow \mathbb{R}$,
we have $\frac{d(\sum_{i=1}^N u_i)}{dt}=0$ along the combinatorial $p$-th Calabi flow (\ref{Eq: PCCF}).
This implies $\sum_{i=1}^Nu_i$ is invariant along the combinatorial $p$-th Calabi flow (\ref{Eq: PCCF}).
\qed

Lemma \ref{Lem: invariant 1} implies that the solutions of the combinatorial Calabi flow (\ref{Eq: CCF}), the $2s$-th order fractional combinatorial Calabi flow (\ref{Eq: FCCF}) and  the combinatorial $p$-th Calabi flow (\ref{Eq: PCCF}) stay in the hyperplane $\Sigma_0:=\{u\in \mathbb{R}^N|\sum^N_{i=1}u_i=\sum^N_{i=1}u_i(0)\}$.
\\

\noindent\textbf{Proof of Theorem \ref{Thm: main 1}:}
Suppose the solution $u(t)$ of the combinatorial Calabi flow (\ref{Eq: CCF}) converges to $\overline{u}$ as $t\rightarrow +\infty$, then
$K(\overline{u})=\lim_{t\rightarrow +\infty}K(u(t))$ by the $C^1$-smoothness of $K$.
Furthermore, there exists a sequence $t_n\in(n,n+1)$ such that for any $i\in V$,
\begin{equation*}
u_i(n+1)-u_i(n)=u'_i(t_n)=\Delta_\mathcal{T}(K(u(t_n))
-\overline{K})_i\rightarrow 0,\ \text{as}\ n\rightarrow +\infty.
\end{equation*}
This implies that $K(\overline{u})-\overline{K}=\lim_{n\rightarrow +\infty}(K(u(t_n))-\overline{K})$
is in the kernel of the discrete Laplace operator $\Delta_{\mathcal{T}}$.
Therefore, by Lemma \ref{Lem: Xu}, we have $K(\overline{u})-\overline{K}=c\mathbf{1}^\mathrm{T}$ for some $c\in \mathbb{R}$.
Note that $\sum_{i=1}^N(K_i(\overline{u})-\overline{K}_i)
=2\pi\chi(S)-2\pi\chi(S)=0$.
This implies $K(\overline{u})-\overline{K}=0$ and $\overline{u}$ is a discrete conformal factor with the combinatorial curvature $\overline{K}$.

Similarly, if the solution of the $2s$-th order fractional combinatorial Calabi flow (\ref{Eq: FCCF}) converges,
then
\begin{equation*}
u_i(n+1)-u_i(n)=u'_i(t_n)
=\Delta^s_\mathcal{T}(K(u(t_n))
-\overline{K})_i\rightarrow 0,\ \text{as}\ n\rightarrow +\infty.
\end{equation*}
By the definition of the $2s$-th order fractional discrete Laplace operator in (\ref{Eq: fractional Laplace}),
$\Delta^s_{\mathcal{T}}$ is symmetric and negative semi-definite with rank $N-1$ and kernel $\{c\mathbf{1}^\mathrm{T}\in \mathbb{R}^N|c\in \mathbb{R}\}$ on $\Omega^\mathcal{T}$.
This implies that $K(\overline{u})-\overline{K}=\lim_{n\rightarrow +\infty}(K(u(t_n))-\overline{K})=c\mathbf{1}^\mathrm{T}$ for some $c\in \mathbb{R}$.
The rest of proof is paralleling to the case of the combinatorial Calabi flow, so we omit it.

Suppose there exists a discrete conformal factor $\overline{u}$ with the combinatorial curvature $\overline{K}$.
For the combinatorial Calabi flow (\ref{Eq: CCF}),
set $\Gamma(u)=\Delta_{\mathcal{T}}(K-\overline{K})$.
Then $D\Gamma|_{u=\overline{u}}=-L^2$ is negative semi-definite with kernel $\{c\mathbf{1}^\mathrm{T}\in \mathbb{R}^N|c\in \mathbb{R}\}$ by Lemma \ref{Lem: Xu}. Note that the kernel is perpendicular to the hyperplane $\Sigma_0$.
By Lemma \ref{Lem: invariant 1}, this implies that $\overline{u}$ is a local attractor of (\ref{Eq: CCF}).
Then the conclusion follows from Lyapunov Stability Theorem (\cite{Pontryagin}, Chapter 5).

Similarly, for the $2s$-th order fractional combinatorial Calabi flow (\ref{Eq: FCCF}),
set $\Gamma(u)=\Delta_\mathcal{T}^s(K-\overline{K})$.
Then $D\Gamma|_{u=\overline{u}}=-L^{s+1}$ restricted to the hypersurface $\Sigma_0$ is negative definite,
which implies that $\overline{u}$ is a local attractor of (\ref{Eq: FCCF}). Then the conclusion follows from Lyapunov Stability Theorem (\cite{Pontryagin}, Chapter 5).
\qed

For general initial decorated PE metrics on $(S,\mathcal{T},I)$,
Theorem \ref{Thm: main 1} does not give the longtime existence and convergence of the solutions of the combinatorial Calabi flow (\ref{Eq: CCF}), the $2s$-th order fractional combinatorial Calabi flow (\ref{Eq: FCCF}) and the combinatorial $p$-th Calabi flow (\ref{Eq: PCCF}).
To study the longtime behavior of these combinatorial curvature flows for general initial decorated PE metrics,
we need to analyze the discrete Laplace operator $\Delta_\mathcal{T}$.
Define the decorated cotan-weights $w: E\rightarrow \mathbb{R}$ by
\begin{equation*}
w_{ij}=\cot \alpha_{ij}^k+\cot \alpha_{ij}^l
=\frac{\sin(\alpha_{ij}^k+\alpha_{ij}^l)}
{\sin\alpha_{ij}^k\sin\alpha_{ij}^l}
\end{equation*}
for any edge $\{ij\}\in E$,
where $\alpha_{ij}^k$ is the interior intersection angle of the face-circle $C_{ijk}$ and the edge $\{ij\}$.
Note that $w_{ij}$ may be negative and unbounded along these combinatorial curvature flows for general initial decorated PE metrics.
Moreover, the weight $w_{ij}\geq0$ is equivalent to $\alpha_{ij}^k+\alpha_{ij}^l\leq \pi$,
which is the weighted Delaunay condition.
Denote $r_{ij}$ as half of the distance of the two intersection points of the face-circle $C_{ijk}$ and the edge $\{ij\}$.
The signed distance $d_{ij}^k$ of the center $c_{ijk}$ of the face-circle $C_{ijk}$ to the edge $\{ij\}$ is given by
$d_{ij}^k=r_{ij}\cot\alpha_{ij}^k$,
which is defined to be positive if the center is on the same side of the line determined by $\{ij\}$ as the triangle $\{ijk\}$ and negative otherwise (or zero if the center is on the line).
Please refer to Figure \ref{figure}.
By direct calculations, we have
\begin{equation}\label{Eq: F1}
\frac{\partial K_i}{\partial u_j}
=-(\frac{d_{ij}^k}{l_{ij}}+\frac{d_{ij}^l}{l_{ij}})
=-(\frac{r_{ij}\cot\alpha_{ij}^k}{l_{ij}}
+\frac{r_{ij}\cot\alpha_{ij}^l}{l_{ij}})
=-\frac{r_{ij}}{l_{ij}}w_{ij}.
\end{equation}
This implies $\frac{\partial K_i}{\partial u_j}\leq0$ under the weighted Delaunay condition.
One can also refer to Subsection 3.3 in \cite{BL} for the equation (\ref{Eq: F1}).
It is natural to require the triangulation $\mathcal{T}$ to be weighted Delaunay,
otherwise the coefficients of the discrete Laplace operator $\Delta_{\mathcal{T}}$ will be negative.
Denote $\mathcal{C}_\mathcal{T}(dist_{S},r)\subseteq \mathbb{R}^N$ as the set of discrete conformal factors $u$ such that $\mathcal{T}$ is still a weighted Delaunay triangulation after discrete conformally changing of the decorated metric by (\ref{Eq: DCE2}), which is a subspace of the admissible space $\Omega^\mathcal{T}$.
Bobenko-Lutz \cite{BL} proved that the space $\mathcal{C}_\mathcal{T}(dist_{S},r)$ is homeomorphic to a polyhedral cone (with its apex removed) and its interior is homeomorphic to $\mathbb{R}^N$.

\begin{remark}\label{Rmk: 2}
Under the conditions in Theorem \ref{Thm: main 1}, if the triangulation $\mathcal{T}$ is weighted Delaunay
along the combinatorial $p$-th Calabi flow (\ref{Eq: PCCF}),
and the solution $u(t)$ of the combinatorial $p$-th Calabi flow (\ref{Eq: PCCF}) converges to $\overline{u}$,
then there exists a discrete conformal factor on $(S, \mathcal{T}, I)$ with the combinatorial curvature $\overline{K}$.
Indeed, there exists a sequence $t_n\in(n,n+1)$ such that
for any $i\in V$,
\begin{equation*}
u_i(n+1)-u_i(n)=u'_i(t_n)
=\Delta_{p,\mathcal{T}}(K(u(t_n))
-\overline{K})_i\rightarrow 0,\ \text{as}\ n\rightarrow +\infty.
\end{equation*}
Set $\widetilde{K}=\lim_{n\rightarrow +\infty}(K(u(\xi_n))-\overline{K})
=K(\overline{u})-\overline{K}$,
then $\Delta_{p,\mathcal{T}}\widetilde{K}=0$.
By the following formula obtained by Lin-Zhang (\cite{L-Z}, Lemma 5.5)
\begin{equation}\label{Eq: LZ}
\sum_{i=1}^Nf_i\Delta_{p,\mathcal{T}}f_i
=\frac{1}{2}\sum_{i=1}^N\sum_{j\sim i}\frac{\partial K_i}{\partial u_j}|f_j-f_i|^p
\end{equation}
for any $f: V\rightarrow \mathbb{R}$,
we have
\begin{equation*}
0=\widetilde{K}^\mathrm{T}\Delta_{p,\mathcal{T}}\widetilde{K}
=\sum_{i}\widetilde{K}_i\Delta_{p,\mathcal{T}}\widetilde{K}_i
=\frac{1}{2}\sum_{i=1}^N\sum_{j\sim i}(\frac{\partial K_i}{\partial u_j})|\widetilde{K}_i
-\widetilde{K}_j|^p.
\end{equation*}
Since $\frac{\partial K_i}{\partial u_j}\leq0$ under the weighted Delaunay condition,
then $\widetilde{K}_i=\widetilde{K}_j$ for the edges $\{ij\}\in E$ with $\frac{\partial K_i}{\partial u_j}<0$.
Note that the edges with $\frac{\partial K_i}{\partial u_j}<0$ correspond to the edges of the canonical weighted Delaunay tessellation
of the PE surface determined by $\overline{u}$.
They form a connected graph connecting all the vertices in $V$.
Therefore, we still have $\widetilde{K}\equiv c$ for some constant $c\in \mathbb{R}$,
which implies $K(\overline{u})-\overline{K}=c\mathbf{1}^\mathrm{T}$.
The rest of proof is similar that of Theorem \ref{Thm: main 1}, so we omit it.
\end{remark}

\section{Combinatorial curvature flows for variable triangulations}\label{Sec: variable triangulation}

\subsection{Bobenko-Lutz's work on discrete uniformization theorem for decorated PE metrics}\label{subsec: DCT}
To analyze the longtime behavior of the combinatorial Calabi flow with surgery, the $2s$-th order fractional combinatorial Calabi flow with surgery and the combinatorial $p$-th Calabi flow with surgery,
we need the discrete conformal theory for decorated PE metrics recently established by Bobenko-Lutz \cite{BL}.
The discrete conformal theory established by Bobenko-Lutz \cite{BL} also applies to Luo's vertex scalings,
and hence generalizes the discrete conformal theory established by Gu-Luo-Sun-Wu \cite{Gu1}.
In this subsection, we briefly recall Bobenko-Lutz's discrete conformal theory related to the inversive distance circle packings. Please refer to \cite{BL} for more details.

\begin{definition}[\cite{BL}, Definition 4.11]\label{Def: GDCE}
Two decorated PE metrics $(dist_{S},r)$ and $(\widetilde{dist}_{S},\widetilde{r})$ on the marked surface $(S,V)$ are discrete conformal equivalent if and only if there is a sequence of triangulated decorated PE surfaces
$(\mathcal{T}^0,l^0,r^0),...,(\mathcal{T}^N,l^N,r^N)$ such that
\begin{description}
  \item[(i)] the decorated PE metric of $(\mathcal{T}^0,l^0,r^0)$ is $(dist_{S},r)$ and the decorated PE metric of $(\mathcal{T}^N,l^N,r^N)$ is $(\widetilde{dist}_{S},\widetilde{r})$,
  \item[(ii)] each $\mathcal{T}^n$ is a weighted Delaunay triangulation of the decorated PE surface $(\mathcal{T}^n,l^n,r^n)$,
  \item[(iii)] if $\mathcal{T}^n=\mathcal{T}^{n+1}$, then there is a discrete conformal factor $u\in \mathbb{R}^N$ such that $(\mathcal{T}^n,l^n,r^n)$ and $(\mathcal{T}^{n+1},l^{n+1},r^{n+1})$ is related by (\ref{Eq: DCE1}) and (\ref{Eq: DCE2}),
  \item[(iv)] if $\mathcal{T}^n\neq\mathcal{T}^{n+1}$, then $\mathcal{T}^n$ and $\mathcal{T}^{n+1}$ are two different weighted Delaunay triangulations of the same decorated PE surface.
\end{description}
\end{definition}

Using the new definition of discrete conformality, Bobwnko-Lutz \cite{BL} proved the following theorem for decorated PE metrics.

\begin{theorem}[\cite{BL}, Theorem A]\label{Thm: BL}
Let $(dist_{S},r)$ be a decorated PE metric on the marked surface $(S,V)$. Then
\begin{description}
  \item[(i)] there exists a decorated PE metric discrete conformal equivalent to $(dist_{S},r)$ realizing $\overline{K}: V\rightarrow (-\infty, 2\pi)$ if and only if $\overline{K}$ satisfies the discrete Gauss-Bonnet foumula (\ref{Eq: Gauss-Bonnet}).
  \item[(ii)] for each $\overline{K}: V\rightarrow (-\infty, 2\pi)$, there exists at most one decorated PE metric discrete conformal equivalent to $(dist_{S},r)$ realizing $\overline{K}$, up to scaling.
\end{description}
\end{theorem}

For any decorated PE surface, there exists a unique complete hyperbolic surface $\Sigma_g$, i.e., the surface induced by any triangular refinement of its unique weighted Delaunay tessellation.
It is homeomorphic to $S\backslash V$ and called the fundamental discrete conformal invariant of the PE metric $(dist_{S},r)$.
The decoration of $\Sigma_g$ is denoted by $\omega:=e^{h}$ and here the height $h$ is related to $u$ by
\begin{equation}\label{Eq: h u}
dh_i=-du_i.
\end{equation}
The canonical weighted Delaunay tessellation of $\Sigma_g$ is denoted by $\mathcal{T}_{\Sigma_g}^\omega$.
Bobenko-Lutz \cite{BL} defined the following set
\begin{equation*}
\mathcal{D}_\mathcal{T}(\Sigma_g)
=\{\omega\in \mathbb{R}_{>0}^N|\mathcal{T}\ \text{refines}\ \mathcal{T}_{\Sigma_g}^\omega\}
\end{equation*}
and proved the following proposition.
\begin{proposition}[\cite{BL}, Proposition 4.3]\label{Prop: finite decomposition}
Given a complete hyperbolic surface with ends $\Sigma_g$.
There is only a finite number of geodesic tessellations $\mathcal{T}_1,...,\mathcal{T}_m$ of $\Sigma_g$ such that $\mathcal{D}_{\mathcal{T}_n}(\Sigma_g)$ $(n=1,...,m)$ is non-empty.
In particular, $\mathbb{R}^N_{>0}
=\bigcup_{n=1}^m\mathcal{D}_{\mathcal{T}_n}(\Sigma_g)$.
\end{proposition}

The set of all heights $h$ of convex polyhedral cusps over the triangulated hyperbolic surface $(\Sigma_g,\mathcal{T})$ is denoted by $\mathcal{P}_\mathcal{T}(\Sigma_g)\subseteq \mathbb{R}^N$.

\begin{proposition}[\cite{BL}, Proposition 4.9]\label{Prop: some spaces}
Given a decorated PE-metric $(dist_{S},r)$ on the marked surface $(S,V)$.
Then $\mathcal{C}_\mathcal{T}(dist_{S},r)$, $\mathcal{P}_\mathcal{T}(\Sigma_g)$ and $\mathcal{D}_\mathcal{T}(\Sigma_g)$ are homeomorphic.
\end{proposition}

Therefore, $\mathbb{R}^N_{>0}
=\bigcup_{n=1}^m\mathcal{C}_{\mathcal{T}_n}(dist_{S},r)$.
A geodesic triangulation $\mathcal{T}$, which refines the unique canonical weighted Delaunay tessellation $\mathcal{T}_{\Sigma_g}^\omega$, is called a canonical weighted Delaunay triangulation of $\Sigma_g$ with respect to the weights $\omega$.
One can extend the combinatorial map $K$ to
\begin{equation}\label{Eq: extended K}
\begin{aligned}
\mathbf{K}: \mathbb{R}^N&\rightarrow \mathbb{R}_{>0}^N\\
u&\mapsto K(u),
\end{aligned}
\end{equation}
which is independent of the choice of the canonical weighted Delaunay triangulations.
Let $\Sigma_g$ be a hyperbolic surface and $\overline{K}: V\rightarrow (-\infty, 2\pi)$ satisfies the discrete Gauss-Bonnet formula (\ref{Eq: Gauss-Bonnet}).
Bobenko-Lutz \cite{BL} defined the following discrete Hilbert-Einstein functional (dHE-functional)
\begin{equation}\label{Eq: functional H}
\mathcal{H}_{\Sigma_g,\overline{K}}(h)
=\mathcal{H}_{\Sigma_g,\overline{K},\mathcal{T}}(h)
=-2\mathrm{Vol}(P_h)+\sum_{i\in V}(\Theta_i-\theta_i)h_i
+\sum_{\{ij\}\in E(\mathcal{T})}(\pi-\alpha_{ij})\lambda_{ij}
\end{equation}
on $\mathbb{R}^N$,
where $P_h$ is the convex polyhedral cusp defined by the heights $h\in \mathbb{R}^N$, $\mathrm{Vol}(P_h)$ is the volume of $P_h$, $\mathcal{T}$ is a canonical weighted Delaunay triangulation corresponding to the weights $e^{h_i}$ and $\alpha_{ij}=\alpha_{ij}^k+\alpha_{ij}^l$.
Note that $\Theta_i=2\pi-\overline{K}_i$ and $\theta_i=2\pi-\mathbf{K}_i$.
The dHE-functional $\mathcal{H}_{\Sigma_g,\overline{K}}$ has the following properties.

\begin{proposition}[\cite{BL}, Proposition 4.13]\label{Prop: dHE-functional}
Let $\Sigma_g$ be a hyperbolic surface and $\overline{K}: V\rightarrow (-\infty, 2\pi)$ satisfy the discrete Gauss-Bonnet formula (\ref{Eq: Gauss-Bonnet}).
\begin{description}
  \item[(i)] The dHE-functional $\mathcal{H}_{\Sigma_g,\overline{K}}$ is concave, twice continuously differentiable over $\mathbb{R}^N$ and analytic in each $\mathcal{P}_\mathcal{T}(\Sigma_g)$.
  \item[(ii)] The first derivative of the dHE-functional $\mathcal{H}_{\Sigma_g,\overline{K}}$ is given by
      \begin{equation}\label{Eq: first derivative}
      d\mathcal{H}_{\Sigma_g,\overline{K}}
      =\sum_{i=1}^N(\Theta_i-\theta_i)dh_i
      =\sum_{i=1}^N(\mathbf{K}_i-\overline{K}_i)dh_i.
      \end{equation}
  \item[(iii)] The dHE-functional $\mathcal{H}_{\Sigma_g,\overline{K}}$ is shift-invariant, i.e., for any $t\in \mathbb{R}$,
  \begin{equation*}
  \mathcal{H}_{\Sigma_g,\overline{K}}(h+t\mathbf{1})
  =\mathcal{H}_{\Sigma_g,\overline{K}}(h).
  \end{equation*}
Furthermore, the restriction of $\mathcal{H}_{\Sigma_g,\overline{K}}$ to $\{h\in\mathbb{R}^N|\sum_i^Nh_i=0\}$ is strictly concave and coercive, i.e.,
\begin{equation*}
\lim_{||h||\rightarrow +\infty}\mathcal{H}_{\Sigma_g,\overline{K}}(h)=-\infty.
\end{equation*}
\end{description}
\end{proposition}

By (\ref{Eq: h u}) and (\ref{Eq: first derivative}),
the dHE-functional $\mathcal{H}_{\Sigma_g,\overline{K}}$ is equivalent to the following function up to a constant
\begin{equation}\label{Eq: function F}
F(u)=\int^u\sum_{i=1}^N(\overline{K}_i-\mathbf{K}_i)du_i,
\end{equation}
which is well-defined by Proposition \ref{Prop: some spaces}.

Proposition \ref{Prop: dHE-functional} implies that $\mathbf{K}$ defined on $\mathbb{R}^N$ is a $C^1$-extension of the combinatorial curvature $K$ defined on the space of discrete conformal factors $\mathcal{C}_\mathcal{T}(dist_{S},r)$.
In general, for a marked surface with a decorated PE metric,
the discrete Lapalce operator $\Delta_{\mathcal{T}}$ depends on the geometric triangulation $\mathcal{T}$ of the PE surface.
However, by the fact that different weighted Delaunay triangulations for the same decorated PE metric on a marked surface $(S,V)$ correspond to the same canonical weighted Delaunay tessellation \cite{BL},
if the triangulation $\mathcal{T}$ is weighted Delaunay,
even though there exists different weighted  Delaunay triangulations for the same decorated PE metric on a marked surface $(S,V)$,
the discrete Laplace operator $\Delta_{\mathcal{T}}$ is independent of the weighted Delaunay triangulations of $(l,r)$ on $(S,V)$.
In this sense, the discrete Laplace operator $\Delta_{\mathcal{T}}$ is intrinsic.
Then the discrete Laplace operator $\Delta_\mathcal{T}$ could be extend to the following operator $\Delta$ defined on $\mathbb{R}^N$, which is continuous and piecewise smooth on $\mathbb{R}^N$ as a matrix-valued function of $u$.

\begin{definition}\label{Def: Delta}
Suppose $(S,V)$ is a marked surface with a decorated PE metric $(dist_g,r)$.
The discrete Laplace operator $\Delta$ is defined to be the map
\begin{equation*}
\begin{aligned}
\Delta :\  &\mathbb{R}^N\rightarrow \mathbb{R}^N\\
&f\mapsto \Delta f,
\end{aligned}
\end{equation*}
with
\begin{equation*}
\Delta f_i
=-\sum_{j\in V}\frac{\partial \mathbf{K}_i}{\partial u_j}f_j
=-(\widetilde{L}f)_i
\end{equation*}
for any $f : V \rightarrow \mathbb{R}$, where $\widetilde{L}_{ij}=\frac{\partial \mathbf{K}_i}{\partial u_j}$ is an extension of $L_{ij}=\frac{\partial K_i}{\partial u_j}$.
\end{definition}
Similarly,
for any $s\in \mathbb{R}$, the extended $2s$-th order fractional discrete Laplace operator $\Delta^s$ is defined to be
\begin{equation*}
\Delta^s=-\widetilde{L}^s.
\end{equation*}
For any $p>1$,
the extended discrete $p$-th Laplace operator $\Delta_{p}$ is defined by
\begin{equation*}
\Delta_{p}f_i
=\sum_{j\sim i}(-\frac{\partial \mathbf{K}_i}{\partial u_j})|f_j-f_i|^{p-2}(f_j-f_i).
\end{equation*}
for $f: V\rightarrow \mathbb{R}$.

\subsection{Combinatorial curvature flows with surgery}

\begin{definition}
Suppose $(S,V)$ is a marked surface with a decorated PE metric $(dist_g,r)$.
Let $\overline{K}: V\rightarrow (-\infty, 2\pi)$ be a given function defined on $V$.
The combinatorial Calabi flow with surgery is defined to be
\begin{eqnarray}\label{Eq: CCFS}
\begin{cases}
\frac{du_i}{dt}=\Delta(\mathbf{K}-\overline{K})_i,\\
u_i(0)=u_0.
\end{cases}
\end{eqnarray}
For any $s\in \mathbb{R}$, the $2s$-th order fractional combinatorial Calabi flow with surgery is defined to be
\begin{eqnarray}\label{Eq: FCCFS}
\begin{cases}
\frac{du_i}{dt}=\Delta^s(\mathbf{K}-\overline{K})_i,\\
u_i(0)=u_0.
\end{cases}
\end{eqnarray}
For any $p>1$, the combinatorial $p$-th Calabi flow with surgery is defined to be
\begin{eqnarray}\label{Eq: PCCFS}
\begin{cases}
\frac{du_i}{dt}=\Delta_{p}(\mathbf{K}-\overline{K})_i,\\
u_i(0)=u_0,
\end{cases}
\end{eqnarray}
\end{definition}

Similar to Lemma \ref{Lem: invariant 1}, we have the following result.
Here, we omit the proof.

\begin{lemma}\label{Lem: invariant 2}
Suppose $(S,V)$ is a marked surface with a decorated PE metric $(dist_g,r)$.
Let $\overline{K}: V\rightarrow (-\infty, 2\pi)$ be a given function defined on $V$ satisfying the discrete Gauss-Bonnet formula (\ref{Eq: Gauss-Bonnet}),
then $\sum_{i=1}^Nu_i(t)$ is invariant along the combinatorial Calabi flow with surgery (\ref{Eq: CCFS}), the $2s$-th order fractional combinatorial Calabi flow with surgery (\ref{Eq: FCCFS}) and the combinatorial $p$-th Calabi flow with surgery (\ref{Eq: PCCFS}).
\end{lemma}

Lemma \ref{Lem: invariant 2} implies that the solutions of the combinatorial Calabi flow with surgery (\ref{Eq: CCFS}), the $2s$-th order fractional combinatorial Calabi flow with surgery (\ref{Eq: FCCFS}) and the combinatorial $p$-th Calabi flow with surgery (\ref{Eq: PCCFS}) stay in the hyperplane $\Sigma_0=\{u\in \mathbb{R}^N|\sum^N_{i=1}u_i=\sum^N_{i=1}u_i(0)\}$.
\\

\noindent\textbf{Proof of Theorem \ref{Thm: main 2}:}
As the function $\overline{K}: V\rightarrow (-\infty, 2\pi)$ satisfies the discrete Gauss-Bonnet formula (\ref{Eq: Gauss-Bonnet}), by Theorem \ref{Thm: BL}, there exists a unique decorated PE metric with the combinatorial curvature $\overline{K}$ and a unique $\overline{u}\in \Sigma_0$ such that $\mathbf{K}(\overline{u})=\overline{K}$.
Set
\begin{equation}\label{Eq: W}
W(u)=-F(u)=\int^u_{\overline{u}}
\sum_{i=1}^N(\mathbf{K}_i-\overline{K}_i)du_i.
\end{equation}
By Proposition \ref{Prop: dHE-functional}, $W$ is a $C^2$-smooth convex function defined on $\mathbb{R}^N$.
Furthermore, $W(\overline{u})=0,\ \nabla W(\overline{u})=0$, $\mathrm{Hess} W\geq0$ and the kernel of
$\mathrm{Hess} W$ is orthogonal to $\Sigma_0$.
This implies $\lim_{||u||\rightarrow +\infty}W(u)|_{\Sigma_0}=+\infty$.
Hence, $W(u)|_{\Sigma_0}$ is proper and $0=W(\overline{u})\leq W(u)$.

\noindent\textbf{(i):}
Along the combinatorial Calabi flow with surgery (\ref{Eq: CCFS}), by Lemma \ref{Lem: Xu}, we have
\begin{equation*}
\frac{dW(u(t))}{dt}
=\sum^N_{i=1}\frac{\partial W}{\partial u_i}\frac{du_i}{dt}
=\sum^N_{i=1}(\mathbf{K}-\overline{K})_i
\Delta(\mathbf{K}-\overline{K})_i
=-(\mathbf{K}-\overline{K})^\mathrm{T}\cdot \widetilde{L}\cdot (\mathbf{K}-\overline{K})\leq0.
\end{equation*}
This implies $0\leq W(u(t))\leq W(u(0))$.
Combining Lemma \ref{Lem: invariant 2} and the properness of $W(u)|_{\Sigma_0}$,
the solution $\{u(t)\}$ of the combinatorial Calabi flow with surgery (\ref{Eq: CCFS}) lies in a compact subset of $\Sigma_0$,
which implies the solution of the combinatorial Calabi flow with surgery (\ref{Eq: CCFS}) exists for all time and $W(u(t))|_{\Sigma_0}$ converges.

Moreover, there exists a sequence $t_n\in(n,n+1)$ such that as $n\rightarrow +\infty$,
\begin{equation*}
\begin{aligned}
&W(u(n+1))-W(u(n))
=(W(u(t))'|_{t_n}
=\nabla W\cdot\frac{du_i}{dt}|_{t_n}\\
=&\sum^N_{i=1}(\mathbf{K}-\overline{K})_i
\Delta(\mathbf{K}-\overline{K})_i|_{t_n}
=-(\mathbf{K}-\overline{K})^\mathrm{T}\cdot \widetilde{L}\cdot (\mathbf{K}-\overline{K})|_{t_n}\rightarrow 0.
\end{aligned}
\end{equation*}
This implies that $\mathbf{K}(\overline{u})-\overline{K}=\lim_{n\rightarrow +\infty}(\mathbf{K}(u(t_n))-\overline{K})$
is in the kernel of the discrete Laplace operator $\Delta$.
Therefore, by Lemma \ref{Lem: Xu}, we have $\mathbf{K}(\overline{u})-\overline{K}=c\mathbf{1}^\mathrm{T}$ for some $c\in \mathbb{R}$.
Note that $\sum_{i=1}^N(\mathbf{K}_i(\overline{u})-\overline{K}_i)
=2\pi\chi(S)-2\pi\chi(S)=0$.
This implies $\mathbf{K}(\overline{u})-\overline{K}=0$.
By $\{u(t)\}\subset\subset \Sigma_0$, there exists $u^*\in \mathbb{R}^N$ and a convergent subsequence of $\{u(t_n)\}$, still denoted as $\{u(t_n)\}$ for simplicity, such that $\lim_{n\rightarrow \infty}u(t_n)=u^*$.
This implies
$\mathbf{K}(u^*)=\lim_{n\rightarrow +\infty}\mathbf{K}(u(t_n))
=\mathbf{K}(\overline{u})$.
Then $u^*=\overline{u}$ by Theorem \ref{Thm: BL}.
Therefore, $\lim_{n\rightarrow \infty}u(t_n)=\overline{u}$.

Set $\Gamma(u)=\Delta (\mathbf{K}-\overline{K})$,
then $D\Gamma|_{u=\overline{u}}$ restricted to the hypersurface $\Sigma_0$ is negative definite,
which implies that $\overline{u}$ is a local attractor of (\ref{Eq: CCFS}).
Then the conclusion follows from Lyapunov Stability Theorem (\cite{Pontryagin}, Chapter 5).

\noindent\textbf{(ii):}
Along the $2s$-th order fractional combinatorial Calabi flow with surgery (\ref{Eq: FCCFS}), by Lemma \ref{Lem: Xu},
we have
\begin{equation*}
\frac{dW(u(t))}{dt}
=\sum^N_{i=1}\frac{\partial W}{\partial u_i}\frac{du_i}{dt}
=\sum^N_{i=1}(\mathbf{K}-\overline{K})_i
\Delta^s(\mathbf{K}-\overline{K})_i
=-(\mathbf{K}-\overline{K})^\mathrm{T}\cdot \widetilde{L}^s\cdot (\mathbf{K}-\overline{K})\leq0.
\end{equation*}
This implies $0\leq W(u(t))\leq W(u(0))$.
By the properness of $W(u)|_{\Sigma_0}$,
the solution $\{u(t)\}$ of the $2s$-th order fractional combinatorial Calabi flow with surgery (\ref{Eq: FCCFS}) lies in a compact subset of $\Sigma_0$.
This implies the solution of the fractional combinatorial Calabi flow with surgery (\ref{Eq: FCCFS}) exists for all time.
By Lemma \ref{Lem: Xu}, the matrix $\widetilde{L}^{s+1}$ is strictly positive definite on $\Sigma_0$.
By the continuity of the eigenvalues of $\widetilde{L}^{s+1}$,
there exists $\lambda_0>0$ such that the non-zero eigenvalues $\lambda$ of $\widetilde{L}^{s+1}$ satisfy $\lambda>\lambda_0$ along the
$2s$-th order fractional combinatorial Calabi flow with surgery (\ref{Eq: FCCFS}).
Therefore, for the combinatorial Calabi energy
$\mathcal{C}(t)=||\mathbf{K}-\overline{K}||^2$, we have
\begin{equation*}
\frac{d\mathcal{C}(u(t))}{dt}
=\sum_{i=1}^N\frac{\partial\mathcal{C}}{\partial u_i}\frac{d u_i}{dt}
=-2(\mathbf{K}-\overline{K})^\mathrm{T}\cdot \widetilde{L}^{s+1}\cdot(\mathbf{K}-\overline{K})\leq -2\lambda_0\mathcal{C}(u(t)),
\end{equation*}
which implies
$\mathcal{C}(u(t))\leq e^{-2\lambda_0t}\mathcal{C}(0)$.
As $\mathbf{K}|_{\Sigma_0}$ is a $C^1$-diffeomorphism from $\Sigma_0$ to $\mathbf{K}(\Sigma_0)$ by (\ref{Eq: extended K}) and Lemma \ref{Lem: Xu},
we have
\begin{equation*}
||u(t)-\overline{u}||^2
\leq C_1||\mathbf{K}(t)-\overline{K}||^2
\leq C_1e^{-2\lambda_0t}||\mathbf{K}(0)-\overline{K}||^2
\leq C_2e^{-2\lambda_0t}
\end{equation*}
for some positive constants $C_1,C_2$.

\noindent\textbf{(iii):}
The equation (\ref{Eq: F1}) implies $\frac{\partial K_i}{\partial u_j}
=-\frac{r_{ij}}{l_{ij}}w_{ij}\leq0$ under the weighted Delaunay condition.
Then $\widetilde{L}_{ij}=\frac{\partial \mathbf{K}_i}{\partial u_j}\leq0$ by Definition \ref{Def: Delta}.
Along the combinatorial $p$-th Calabi flow with surgery (\ref{Eq: PCCFS}), we have
\begin{equation*}
\begin{aligned}
\frac{dW(u(t))}{dt}
=&\sum^N_{i=1}\frac{\partial W}{\partial u_i}\frac{du_i}{dt}\\
=&\sum^N_{i=1}(\mathbf{K}-\overline{K})_i
\Delta_p(\mathbf{K}-\overline{K})_i\\
=&\frac{1}{2}\sum_{i=1}^N\sum_{j\sim i}(\frac{\partial \mathbf{K}_i}{\partial u_j})|(\mathbf{K}-\overline{K})_i
-(\mathbf{K}-\overline{K})_j|^p\\
\leq&0,
\end{aligned}
\end{equation*}
where (\ref{Eq: LZ}) is used in the third line.
This implies $0\leq W(u(t))\leq W(u(0))$.
By the properness of $W(u)|_{\Sigma_0}$,
the solution $\{u(t)\}$ of the combinatorial $p$-th Calabi flow with surgery (\ref{Eq: PCCFS}) lies in a compact subset of $\Sigma_0$.
This implies the solution of the combinatorial $p$-th Calabi flow with surgery (\ref{Eq: PCCFS}) exists for all time and $W(u(t))$ converges.

Moreover, there exists a sequence $t_n\in(n,n+1)$ such that as $n\rightarrow +\infty$,
\begin{equation*}
\begin{aligned}
&W(u(n+1))-W(u(n))
=\nabla W\cdot\frac{du_i}{dt}|_{t_n}
=\sum^N_{i=1}(\mathbf{K}-\overline{K})_i
\Delta_p(\mathbf{K}-\overline{K})_i|_{t_n}\\
=&\frac{1}{2}\sum_{i=1}^N\sum_{j\sim i}(\frac{\partial \mathbf{K}_i}{\partial u_j})|(\mathbf{K}-\overline{K})_i
-(\mathbf{K}-\overline{K})_j|^p|_{t_n}\rightarrow 0.
\end{aligned}
\end{equation*}
By the similar arguments in Remark \ref{Rmk: 2},
we have $\mathbf{K}(u(t_n))-\overline{K}=c\mathbf{1}^\mathrm{T}$ for some constant $c\in \mathbb{R}$.
Since $\overline{K}$ satisfies the discrete Gauss-Bonnet formula (\ref{Eq: Gauss-Bonnet}),
then $\lim_{n\rightarrow +\infty}\mathbf{K}_i(u(t_n))=\overline{K}_i
=\mathbf{K}_i(\overline{u})$ for all $i\in V$.
By $\{u(t)\}\subset\subset \Sigma_0$,
there exists $u^*\in \mathbb{R}^N$ and a convergent subsequence $\{u(t_{n_k})\}$ of $\{u(t_n)\}$ such that
$\lim_{k\rightarrow \infty}u(t_{n_k})=u^*$.
This implies $\mathbf{K}(u^*)=\lim_{k\rightarrow +\infty}\mathbf{K}(t_{n_k})=\mathbf{K}(\overline{u})$.
Then $u^*=\overline{u}$ by Theorem \ref{Thm: BL}.
Therefore, $\lim_{k\rightarrow\infty}u(t_{n_k})=\overline{u}$.

We use Lin-Zhang's trick in \cite{L-Z} to prove  $\lim_{t\rightarrow\infty}u(t)=\overline{u}$.
Suppose otherwise, there exists $\delta>0$ and $\xi_n\rightarrow +\infty$ such that
$|u(\xi_n)-u^*|>\delta$.
This implies $\{u(\xi_n)\}\subseteq \Sigma_0\backslash B(u^*,\delta)$,
where $B(u^*,\delta)$ is a ball centered at $u^*$ with radius $\delta$.
Hence, for any $u\in \Sigma_0\backslash B(u^*,\delta)$, $W(u)\geq C>0$.
Then $W(u(\xi_n))\geq C>0$.
Since $W(u(t))$ converges and $\lim_{k\rightarrow\infty}u(t_{n_k})=\overline{u}$,
then $W(+\infty)=\lim_{k\rightarrow\infty}W(u(t_{n_k}))
=W(\overline{u})=0$.
Hence, $\lim_{n\rightarrow\infty}W(u(\xi_n))=W(+\infty)=0$.
This is a contradiction.
\qed

\end{document}